\documentclass[journal]{IEEEtran}
\usepackage{amsfonts,amssymb,amsmath}
\usepackage{algorithm,algorithmic}
\usepackage{graphicx}  
\usepackage{epstopdf}
\usepackage{hyperref}
\usepackage{xcolor}
\usepackage{theorem}
\usepackage{cite}
\usepackage{url}
\usepackage{multirow}
\usepackage{booktabs}
\usepackage{multicol}
\usepackage{footmisc}
\usepackage{enumitem}
\usepackage{subfig}
\usepackage{xspace}
\usepackage{bm}
\allowdisplaybreaks[4]
\hypersetup{hidelinks}
\setitemize[1]{itemsep=4.5pt,partopsep=1pt,parsep=\parskip,topsep=4.5pt}

\newtheorem{theorem}{Theorem}
\newtheorem{remark}{Remark}

\newtheorem{definition}{Definition}

\usepackage[acronym]{glossaries}
\newacronym{nlp}{\textsc{nlp}}{nonlinear programming}
\newacronym{qp}{\textsc{qp}}{quadratic program}
\newacronym{ocd}{\textsc{ocd}}{Optimality Condition Decomposition}
\newacronym{app}{\textsc{app}}{Auxiliary Problem Principle}
\newacronym{sqp}{\textsc{sqp}}{Sequential Quadratic Programming}
\newacronym{qcqp}{\textsc{qcqp}}{Quadratically Constrained Quadratic Program} 
\newacronym{rapidpf}{rapid\textsc{pf}}{rapid prototyping for distributed Power Flow}
\newacronym{admm}{\textsc{admm}}{Alternating Direction Method of Multipliers}
\newacronym{aladin}{\textsc{aladin}}{Augmented Lagrangian based Alternating Direction Inexact Newton method}
\newacronym{bfgs}{\textsc{bfgs}}{Broyden–Fletcher–Goldfarb–Shanno}
\newacronym{pcc}{PPC}{point of common coupling}
\newacronym{vsc}{VSC}{voltage source converter}
\newacronym{mtdc}{MTDC}{multiterminal high voltage direct current}
\newacronym{hvdc}{HVDC}{high voltage direct current }
\newacronym{hvac}{HVAC}{high voltage alternating current }
\newacronym{lcc}{LCC}{line commutated converters}
\newacronym{lcc-hvdc}{LCC-HVDC}{\acrfull{lcc} based multiterminal high voltage direct current}
\newacronym{vsc-mtdc}{VSC-MTDC}{\acrfull{vsc} based multiterminal high voltage direct current}
\newacronym{igbt}{IGBT}{insulated gate bipolar transistor}
\newacronym{vsc-hvdc}{VSC-HVDC}{\acrfull{vsc} based high voltage direct current}
\newacronym{opf}{OPF}{optimal power flow}
\newacronym{pu}{p.u.}{per-unit system}
\newacronym{pwm}{PWM}{pulse-width modulation}
\newcommand{\norm}[1]{\left\lVert#1\right\rVert}
\newcommand{\matlab}{\texttt{Matlab}\xspace}
\newcommand{\matpower}{\texttt{Matpower}\xspace}
\newcommand{\casadi}{\textsc{c}as\textsc{ad}i\xspace}
\newcommand{\ipopt }{\texttt{ipopt}\xspace}

\newcommand{\tb}{\textcolor{black}}
\newcommand{\tbb}{\textcolor{black}}

\begin{document}
\title{\LARGE \bf Distributed Optimal Power Flow for VSC-MTDC Meshed AC/DC Grids Using ALADIN}
\author{
Junyi~Zhai, \IEEEmembership{Member, IEEE}, Xinliang~Dai, 
Yuning~Jiang, \IEEEmembership{Member, IEEE}, Ying Xue, \IEEEmembership{Senior Member, IEEE}\\
Veit Hagenmeyer, \IEEEmembership{Member, IEEE}, Colin N. Jones, \IEEEmembership{Member, IEEE},
Xiao-Ping~Zhang,~\IEEEmembership{Fellow,~IEEE}

\thanks{The first two authors contributed equally. This work was supported in part from the Youth Program of Natural Science Foundation of Jiangsu Province (BK20210103) and the Swiss National Science Foundation under the RISK project (Risk Aware Data Driven Demand Response, grant number 200021 175627) (Corresponding author: Yuning Jiang)}
\thanks{J. Zhai is with College of New Energy, China University of Petroleum (East China), Qingdao, China, and with State Grid (suzhou) City \& Energy Research Institute, Suzhou, China (e-mail: {\tt zhaijunyi@163.com})}
\thanks{X. Dai and V. Hagenmeyer are with Institute for Automation and Applied Informatics, Karlsruhe Institute of Technology, Germany. (e-mail: {\tt xingliang.dai, veit hagenmeyer@kit.edu})}
\thanks{Y. Jiang and C. Jones are with Automatic Control Laboratory, EPFL, Switzerland.  (e-mail: {\tt yuning.jiang, colin.jones@epfl.ch})}
\thanks{Y. Xue is with School of Electric Power Engineering, South China University of Technology, Guangzhou, China, and was with Department of Electronic, Electrical and Systems Engineering, University of Birmingham, Birmingham, United Kingdom (e-mail: {\tt dr.yingxue@foxmail.com})}
\thanks{X. Zhang is with Department of Electronic, Electrical and Systems Engineering, University of Birmingham, Birmingham, United Kingdom (e-mail: {\tt x.p.zhang@bham.ac.uk})}}

\maketitle

\begin{abstract}
The increasing application of \acrfull{vsc-hvdc} technology in power grids has raised the importance of incorporating DC grids and converters into the existing transmission network. This poses significant challenges in dealing with the resulting \acrfull{opf} problem. In this paper, a recently proposed nonconvex distributed optimization algorithm---\acrfull{aladin}, is tailored to solve the nonconvex AC/DC \acrshort{opf} problem for emerging \acrfull{vsc-mtdc} meshed AC/DC hybrid systems. The proposed scheme decomposes this AC/DC hybrid \acrshort{opf} problem and handles it in a fully distributed way. Compared to the existing state-of-art \acrfull{admm}, which is in general, not applicable for nonconvex problems, \acrshort{aladin} has a theoretical convergence guarantee. Applying these two approaches to \acrshort{vsc-mtdc} coupled with an IEEE benchmark AC power system illustrates that the tailored \acrshort{aladin} outperforms \acrshort{admm} in convergence speed and numerical robustness. 
\end{abstract}


\IEEEpeerreviewmaketitle

\section{Introduction}
\label{sec::Intro}

Due to AC grid expansion being limited by legislation and the capacity of long-distance transmission, \acrfull{hvdc}---especially \acrfull{vsc-mtdc} technology---is being considered as an alternative solution. \tb{\acrshort{hvdc} applications have traditionally been restricted to the transmission of power between two buses in the AC grids, which are almost exclusively built using thyristor based \acrshort{lcc-hvdc}s. The main drawback of \acrshort{lcc-hvdc} is that it cannot independently control the active and reactive power and suffers from commutation failure under inverter AC faults. Therefore it cannot be connected to very weak AC systems. On the contrary, \acrshort{vsc}s can connect to very weak AC systems, they are able to independently control active and reactive power. This is benefical for controlling the voltage and frequency when connecting to renewable dominated grid. Another great advantage of \acrshort{vsc} over \acrshort{lcc} is that it can be employed in systems with more than two terminals, thereby forming a multi-terminal DC system. The world’s first commercial operational \acrshort{mtdc} system, Italy-Corsica-Sardinia, was built in 1988 \cite{Billon1989}. Now there are 35 \acrshort{vsc-hvdc} systems in operation and 51 projects planned until 2019 \cite{Barnes2019,ZhaiEPSR2021}. In Europe, a pan-European supergrid project is proposed for off-shore wind power transmission utilizing \acrshort{vsc-hvdc} \cite{Rouzbehi2015}. In China, commissioned \acrshort{vsc-mtdc} projects, such as Nanao \acrshort{vsc-mtdc} project \cite{WangYizhen2016}, Zhoushan \acrshort{vsc-mtdc} project \cite{Kirakosyan2018}, and Zhangbei \acrshort{vsc-mtdc} project~\cite{Kong2017} are built. These projects suggest that the \acrshort{vsc-mtdc} system has become a realistic solution for domestic and international power transmission.}

\tb{ Such \acrshort{vsc-mtdc} system provides several advantages, including increased reliability of power transmission (by providing alternative paths), improved balancing services through AC interconnections (by controlling power flow), and reduced generation variability (by sharing the diverse portfolio of intermittent energy resources)~\cite{Zhai2021,MengKe2017,zhouming2018}. Given that \acrshort{mtdc} grids offer unique capability in terms of regulating power flow, which means that the \acrshort{vsc-mtdc} meshed AC/DC system can be operated more flexibly and cost-effectively compared with point-to-point \acrshort{hvdc} connections.} As a result, the \acrfull{opf} problem of meshed AC/DC systems is becoming an urgent task to be investigated carefully. For the purpose of economic efficiency, \acrshort{opf} is used to determine the operating points of the meshed AC/DC grid. In~\cite{1988The}, the first \acrshort{opf} algorithm was introduced for AC grids incorporating point-to-point \acrshort{hvdc} connections, which is formulated as a quadratic programming problem without consideration for converter losses, converter transformers, and filters. In~\cite{Wiget2012} and~\cite{Aragues2012}, a nonlinear AC/DC \acrshort{opf} model was proposed including a quadratic loss model for \acrshort{hvdc} converters.  In fact, the converter losses can add up to a significant fraction of the overall system losses. Thus, the converter losses should be carefully considered in the AC/DC \acrshort{opf} problem. In various literature~\cite{Cao2013,Yang2018,Feng2014,Ergun2019,Khan2019}, the converter losses are more accurately approximated by a quadratic function of its current magnitude. The AC/DC \acrshort{opf} problem is a nonconvex problem due to the nonlinear power flow equations and the highly nonlinear operating constraints imposed by the converters. Several methods have been proposed to address the AC/DC \acrshort{opf} with \acrshort{vsc-hvdc} connections such as heuristic and interior point methods~\cite{Cao2013}, Newton-Raphson method \cite{Pizano2007}, second-order cone programming~\cite{Baradar2013}, sequential quadratic programming methods~\cite{Rabiee2015}, or semidefinite programming relaxation methods~\cite{Bahrami2017}. Although the above literature has investigated the \acrshort{opf} of AC/DC systems, they are formulated in a centralized manner without fully considering the multilevel structure of \acrshort{vsc-mtdc} meshed AC/DC systems.

Although regional grids are physically connected in this context, they make their generating plans independently, with very limited interactions with their neighbours. If a centralized optimization approach is adopted, all the generation information and network topology data needs to be acquired by a single central entity. \tbb{This centralization may create substantial regulatory and political issues because the local system operators have to give up their governance and control to the central entity, which becomes almost impossible under the deregulated electricity market~\cite{Zhai2021}.} Hence, distributed optimization approaches have drawn significant attention in recent times. The most well-known distributed algorithms for tackling the AC \acrshort{opf} problems are \acrfull{ocd}~\cite{Hug2009}, \acrfull{app}~\cite{Baldick1999}, and \acrfull{admm}~\cite{ZhaiIEGS2022,ZhaiITD2022,6748974}. \acrshort{ocd} applies a modified Lagrangian relaxation decomposition. The original problem is partitioned into several subproblems, in which the local variables are decision variables and all foreign variables are fixed to the values of the previous iteration. By penalizing coupling constraints, \acrshort{ocd}~can converge to a solution with moderate accuracy under certain assumptions, e.g., relative weakly coupled subproblems, which cannot be guaranteed in general. 
\tb{
In contrast to \acrshort{ocd}, each subproblem shares the coupling variables with their neighbors and reformulates their objective function by using Augmented Lagrangian Relaxation in the context of both \acrshort{app} and \acrshort{admm}. Compared to~\acrshort{app}, \acrshort{admm} reduces necessary information exchange and only requires neighbor-to-neighbor communication such that \acrshort{admm} surpasses \acrshort{app} in terms of  communication effort. 
In \cite{6748974,Erseghe2015,Guo2017}, \acrshort{admm} has been adopted to deal with the AC \acrshort{opf} problem. However, there are no generic convergence guarantees for AC \acrshort{opf} using the standard \acrshort{admm}. Although~\cite{Erseghe2015} shows the convergence under some technical assumptions, the branch flow limits were not considered and these additional nonlinearities always result in divergence of \acrshort{admm}. Recently, ~\cite{sun2021twolevel} proposed a two-level \acrshort{admm} approach, which as an exception, can deal with AC \acrshort{opf} problem with branch flow limits while guaranteeing convergence. Nevertheless, the AC \acrshort{opf} was formulated as a \acrfull{qcqp} problem at the cost of accuracy, and it numerically converges slowly to a modest accuracy.}

Different from these existing approaches, the \acrfull{aladin} presented recently in~\cite{Boris2016} for distributed nonconvex optimization can provide a local convergence guarantee. As a distributed approach, the local agents' associated subproblems are solved in the parallelizable step of \acrshort{aladin} while an equality constrained \acrfull{qp} is solved in the consensus step, which is in principle close to a \acrfull{sqp} step such that the convergence is sped up. At the expense of communication effort, \acrshort{aladin} obtains locally quadratic convergence for nonconvex problems if suitable Hessian approximations are used, and a global convergence guarantee can also be achieved if the globalization strategy proposed in~\cite[Algorithm~3]{Boris2016} is implemented. 

\tbb{Based on the existing literatures, there are still research gaps that need to be studied.}
\begin{enumerate}
		\item The traditional highly nonlinear and nonconvex \acrshort{opf} problem for meshed AC/DC hybrid systems is usually solved in a centralized manner. In reality, the different regional AC grids and the \acrshort{mtdc} grid are owned and operated by different utilities under the deregulated electricity market, which is an independent decision-making process. Thus, the distributed optimization is much more preferable. Regarding the special structure of \acrshort{vsc-mtdc} meshed AC/DC hybrid systems, it can be partitioned according to the grid-type naturally, and the resulting separated AC/DC \acrshort{opf} problem can be solved by applying a distributed algorithm. 

      \item In a practical \acrshort{vsc-mtdc} meshed AC/DC hybrid systems, considering the AC branch flow limits, the DC branch flow limits, and the dynamic \acrshort{vsc} operating constraints for meshed AC/DC hybrid systems leads to additional highly nonlinear and nonconvex inequality constraints in the resulting optimization problem such that the existing state-of-art distributed algorithms, such as \acrshort{admm}, cannot guarantee the convergence. Recently, an exceptional distributed nonconvex optimization algorithm \acrshort{aladin} has been investigated to deal with distributed AC OPF in~\cite{Engelmann2019} and distributed AC/DC hybrid OPF in~\cite{Meyer2019} with convergence guarantee. The first is our previous work~\cite{Engelmann2019} that studied AC-\acrshort{opf} problem for the AC transmission system, which is the first journal paper applying \acrshort{aladin} to deal with AC-\acrshort{opf} problem in a distributed manner. However, the most challenging AC branch flow limits were not considered. The second is literature~\cite{Meyer2019} that adopted \acrshort{aladin} to deal with the distributed AC/DC \acrshort{opf} problem for an oversimplified AC/DC hybrid system. However, the necessary and most challenging AC branch flow limits, DC branch flow limits, and  highly nonlinear and nonconvex \acrshort{vsc} operating constraints are all ignored. In terms of physical model, ignoring the AC and DC branch flow limits may lead to line overload; ignoring the essential dynamic operating constraints of \acrshort{vsc} stations is unacceptable in the practical operation of \acrshort{vsc-mtdc} meshed AC/DC grids (e.g. the converter power losses can add up to a significant fraction of the overall system losses, which are essential in the AC/DC \acrshort{opf} problem~\cite{Feng2014,YangZhifang2018,Bahrami2017,Ergun2019}). In terms of mathematics, such an oversimplified scheme is much less practical and always leads to relatively weaker couplings between the AC grids and the DC grid. As a result, it is much easier to be solved by using distributed algorithms while heuristic implementation could simply work without encountering numerical challenges.
\end{enumerate}

To fill these research gaps, this paper investigates the potential of using \acrshort{aladin} for the highly nonlinear and nonconvex \acrshort{opf} problem of \acrshort{vsc-mtdc} meshed AC/DC grids. The contributions are summarized as follows.

\begin{enumerate}
\item We consider complicated \acrshort{vsc-mtdc} meshed AC/DC grids including AC/DC interconnections compared to the traditional point-to-point HVDC link. The full nonlinear and  nonconvex AC/DC \acrshort{opf} model is built in an affinely coupled separable form and thus, is feasible for distributed algorithms. Compared to the existing oversimplified model on \acrshort{aladin}-based distributed AC/DC \acrshort{opf}~\cite{Meyer2019}, this paper simultaneously takes the AC branch flow limits, the DC branch flow limits, and the complex dynamic operating model of VSC stations into account. Mathematically, this leads to additional highly nonlinear and nonconvex inequality constraints in the resulting optimization problem. Then, we tailor \acrshort{aladin} algorithm to tackle this challenging nonconvex \acrshort{opf} problem in a distributed manner. At each iteration, the regional subproblems are solved in parallel while an equality constrained \acrshort{qp} is solved in the consensus step, which is in principle close to a \acrshort{sqp} step such that the convergence is sped up. The computation speed of the tailored~\acrshort{aladin} is numerically demonstrated to be faster than traditional centralized method. This indicates that the proposed approach will be of great potential in practical distributed implementation.

\item A detailed numerical investigation of \acrshort{aladin} for the proposed AC/DC \acrshort{opf} problem is presented. Compared to the existing state-of-art \acrshort{admm}, which is, in general, has no convergence guarantee for nonconvex problems, \acrshort{aladin} has theoretical local convergence guarantee. In practice, although~\acrshort{admm} has been applied to solve nonconvex AC~\acrshort{opf}~\cite{Erseghe2015,Guo2017}, one can always observe its divergence when including the branch flow limits. Our numerical comparison illustrates that~\acrshort{aladin} outperforms \acrshort{admm} in all perspective for the proposed nonconvex AC/DC \acrshort{opf} problem. The improved performance comes at the cost of an increased per-step communication effort, which can be reduced by using quasi-Newton Hessian approximation. To this end, the number of iterations is slightly increased but~\acrshort{aladin} is still faster and more accurate than \acrshort{admm}.
\end{enumerate}
\color{black}
The rest of this paper is organized as follows: Section~\ref{sec::decomposition} presents the non-convex AC/DC \acrshort{opf} model. Section~\ref{sec::singleArea} introduces the distributed optimization formulation. Section~\ref{sec::result} presents numerical results. Section~\ref{sec::conc} concludes this paper.

\section{Problem Formulation}
\label{sec::decomposition}
In this section, a mathematical model of the \acrshort{vsc-mtdc} meshed AC/DC grids, e.g., the Zhangbei four terminal \acrshort{mtdc} grid in China~\cite{Kong2017}, is presented. Then, both centralized and distributed OPF formulations are respectively discussed.

\subsection{System Model of VSC-MTDC Meshed AC/DC Grids}
\begin{figure}[htbp!]
\centering
\includegraphics[width=\linewidth]{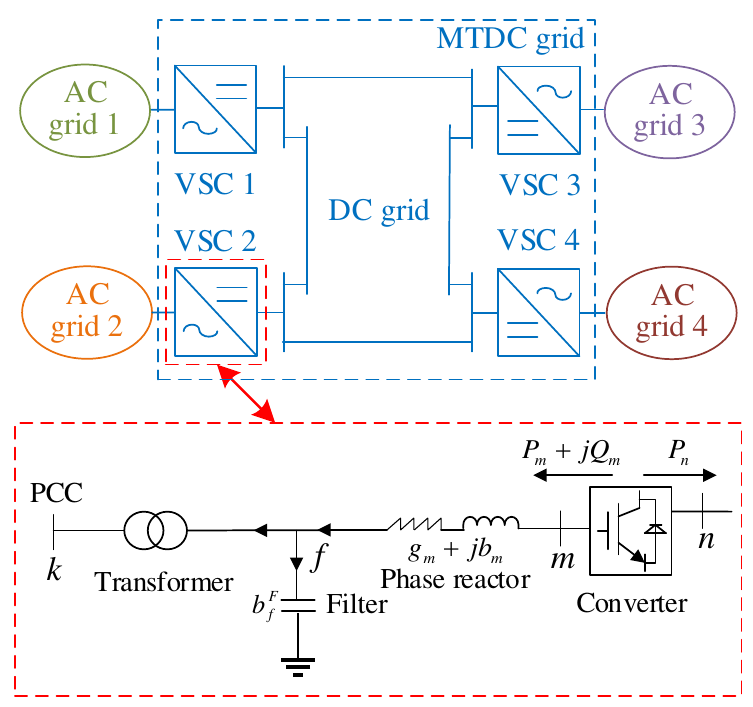}
\caption{VSC station schematic in \acrshort{vsc-mtdc} meshed AC/DC grids with \acrfull{pcc}}
\label{Fg::vsc-mtdc}
\end{figure}

\autoref{Fg::vsc-mtdc} shows the topology of the proposed grid model, \tb{which is abstracted from the Zhangbei four-terminal \acrshort{vsc-mtdc} grid connected with four AC grids~\cite{Kong2017} in Northern China.} The \acrshort{vsc} station is considered with a transformer, AC filter, phase reactor, and converter. Without loss of generality, the \acrshort{vsc} station is assumed to be a two- or three-level converter using the \acrfull{pwm} switching method. We assume there is no generation or load in the \acrshort{mtdc} grid.

Before discussing the grid model further, we introduce some nomenclature. We represent the meshed AC/DC grids by a tuple $(\mathcal{R},\mathcal{N},\mathcal{N}^{ac},\mathcal{N}^{mtdc},\mathcal{L},\mathcal{L}^{ac},\mathcal{L}^{mtdc})$.
Thereby, $\mathcal{R}$ represents the set of all sub-grids, including 4 AC grids and a MTDC grid, i.e., $|\mathcal{R}|=5$,  $\mathcal{N}=\mathcal{N}^{ac}\cup\mathcal{N}^{mtdc}$ denotes the set of all buses, $\mathcal{N}^{ac}$ the set of buses in AC grids, $\mathcal{N}^{mtdc}=\mathcal{N}^{vsc}\cup \mathcal{N}^{dc}$ the set of buses in MTDC grid, $\mathcal{N}^{dc}$ the set of buses in DC grid, $\mathcal{N}^{vsc}=\mathcal{N}_{ac}^{c}\cup \mathcal{N}_{dc}^{c}\cup\mathcal{N}_{ac}^{p}\cup\mathcal{N}_{ac}^{f}$ the set of all VSC buses, $\mathcal{N}_{ac}^{c}$ the set of AC side converter buses, $\mathcal{N}_{dc}^{c}$ the set of DC side converter buses, $\mathcal{N}_{ac}^{p}$ the set of PCC buses, $\mathcal{N}_{ac}^{f}$ the set of filter buses.  $\mathcal{L}=\mathcal{L}^{ac}\cup \mathcal{L}^{mtdc}$ denotes the set of all branches, $\mathcal{L}^{ac}$ the set of branches in AC grids, $\mathcal{L}^{mtdc}$ the set of DC branches in MTDC grid. In \autoref{Fg::vsc-mtdc}, for example, bus $k$ is in set $\mathcal{N}_{ac}^{p}$, bus $f$ is in set $\mathcal{N}_{ac}^{f}$, bus $n$ is in set $\mathcal{N}_{dc}^{c}$, and bus $m$ is in set $\mathcal{N}_{ac}^{c}$.


\subsection{Objective}
The objective of the AC/DC \acrshort{opf} problem is to jointly minimize the total generation cost and the total power losses on the branches and converters, where the total power losses are equal to the total generation minus the total load of the system. The objective of a subproblem can be divided into two parts, the total generation cost given by 
\begin{equation}
\label{eq::C1}
C_1=\sum_{i\in \mathcal{N}^{ac}}\left\{c_{1i}\left( P_{i}^{G} \right)^{2}+c_{2i}P_{i}^{G}+c_{3i}\right\},
\end{equation}
and the total system losses given by
\begin{equation}
\label{eq::C2}
C_2=\sum_{i\in \mathcal{N}^{ac}}\left(P_{i}^{G} -P_{i}^{D}\right),
\end{equation}
where $c_{1i}$, $c_{2i}$, and $c_{3i}$ denote the cost coefficients of generator $i\in\mathcal{N}^{ac}$. $P_i^G$ and $P_i^D$ denote the generator active power output and load at bus $i$, respectively. If bus $i$ is not a generator bus, then, we have $P_i^G=0$ same to the load bus. For this, we only consider the active power losses of the system since the reactive power does not dissipate energy.

The objective function of the AC/DC \acrshort{opf} problem for \acrshort{vsc-mtdc} meshed AC/DC grids is given by
\begin{equation}
\label{eq::C}
C=C_1+\eta C_2,
\end{equation}
where $\eta$ denote a positive scaling coefficient. By increasing $\eta$, the total system losses have a larger weight in the objective function as compared with the generation cost. The typical value of $\eta$ is around the value of coefficients $c_{2i}$, since the total system losses is usually a linear function of the generators’ output. An appropriate value of $\eta$ can enable timely adjustment of control settings to jointly reduce the generation costs, VSC losses and transmission line losses. Thus, it can improve the operational economic efficiency of the entire system.

\subsection{Constraints of AC System}
For the AC/DC \acrshort{opf} problem, the constraints of AC grids consist of power flow balance, branch flow limit and the limit on decision variables, i.e., voltage magnitude of all buses, as well as active and reactive power of generators.
\subsubsection{Nodal power balance of AC grid}  
\begin{subequations}
\label{eq::nodalAC}
\begin{align}
&P_i^G-P_i^D=V_i\sum_{j\in i}\left(G_{ij}\cos \theta_{ij}+B_{ij}\sin \theta_{ij}\right)V_j,\\
&Q_i^G-Q_i^D=V_i\sum_{j\in i}\left(G_{ij}\sin \theta_{ij}-B_{ij}\cos \theta_{ij}\right)V_j,
\end{align}
\end{subequations}
with $i\in \mathcal{N}^{ac}$ and $j\in \mathcal{N}^{ac}\cup\mathcal{N}_{ac}^p$. Notice that the bus $l\in\mathcal{N}_{ac}^p$ belongs to \acrshort{vsc-mtdc} system and branch $(l',l)$ is the linking AC tie-line shown in~\autoref{fig::Fg2}. $G_{ij}$ and $B_{ij}$ denote the real and imaginary part in the admittance matrix. $Q_i^G$ and $Q_i^D$ denote the reactive power output of generator and reactive load at bus $i\in \mathcal{N}^{ac}$. $\theta_{ij}$ denotes the phase angle difference between buses $i$ and $j$. $V_i$ denotes the voltage magnitude of bus $i\in \mathcal{N}^{ac}$.

\subsubsection{Branch flow limit of AC grid}
\begin{subequations}
\label{eq::branAC}
\begin{align}
P_{ij}&=V_i^2g_{ij}-V_i V_j\left(g_{ij}\cos \theta_{ij}+b_{ij}\sin \theta_{ij}\right)\\
Q_{ij}&=-V_i^2b_{ij}+V_i V_j\left(b_{ij}\cos \theta_{ij}-g_{ij}\sin \theta_{ij}\right)\\\label{lineAC}
P_{ij}^2&+Q_{ij}^2 \leq \overline S_{ij}^2&
\end{align}
\end{subequations}
for all $\left(i,j\right) \in \mathcal{L}^{ac}$. Here $g_{ij}$ and $b_{ij}$ denote the conductance and susceptance of branch $\left(i,j\right)\in \mathcal{L}^{ac}$. $P_{ij}$ and $Q_{ij}$ denote the active and reactive power on
branch $\left(i,j\right)\in \mathcal{L}^{ac}$. $\overline S_{ij}$ denotes the maximum allowable apparent power flow through branch $\left(i,j\right)\in \mathcal{L}^{ac}$.

\subsubsection{Limits on voltage magnitude, generator active and reactive power of AC grid} 
\begin{subequations}
\label{eq::actreactAC}
\begin{align}
\underline V_i&\le V_i\le \overline V_i,\quad\;\, i\in\mathcal{N}^{ac},\\
\underline P_i^G&\le P_i^G\le \overline P_i^G,\;\;i\in\mathcal{N}^{ac},\\
\underline Q_i^G&\le Q_i^G\le \overline Q_i^G,\;\;i\in\mathcal{N}^{ac},
\end{align}
\end{subequations}
where $\underline V_i$ and $\overline V_i$ denote the minimum and maximum nodal voltage magnitude at bus $i \in \mathcal{N}^{ac}$. $\underline P_i^G$ and $\overline P_i^G$ denote the minimum and maximum active power of generator $i\in\mathcal{N}^{ac}$. $\underline Q_i^G$ and $\overline Q_i^G$ denote the minimum and maximum reactive power of generator $i\in\mathcal{N}^{ac}$. 

\subsection{Constraints of VSC-MTDC System}
Constraints of the \acrshort{vsc-mtdc} system consist of nodal power balance of \acrshort{mtdc} system, operating limits of converters, as well as limits on branch flow and nodal voltage.
\subsubsection{Nodal power balance inside a VSC station}
The power flow in a \acrshort{vsc} station should satisfy the similar AC power balance at PCC bus $k$ and filter bus $f$, c.f.~\eqref{eq::nodalAC}:
\begin{subequations}
\label{eq::filter}
\begin{align}
0=&V_i\sum_{j\in i}\left(G_{ij}\cos \theta_{ij}+B_{ij}\sin \theta_{ij}\right)V_j,\\
-Q_i^D=&V_i\sum_{j\in i}\left(G_{ij}\sin \theta_{ij}-B_{ij}\cos \theta_{ij}\right)V_j
\end{align}
\end{subequations}
with $i\in\{k,f\}$ and $j\in \{k',k,f,m\}$. Notice that the bus $k'$ belongs to AC system and branch $(k',k)$ is the linking AC tie-line shown in~\autoref{fig::Fg2}. The coupling between AC system and \acrshort{vsc-mtdc} system will be discussed in the rest part of this section.

Regarding the reactive power of the AC filter, the demand at the AC filter bus, i.e., bus $f$, can be expressed as
\begin{equation}
\label{eq::filtAC}
Q_f^D=-V_f^2 b_f^F,\; f\in\mathcal{N}_{ac}^{f}
\end{equation}
where $b_f^F$ denotes the shunt susceptance of the AC filter connected to bus $f$, while there is no reactive power demand at the PCC bus $k$.

\subsubsection{Voltage coupling between the converter AC bus and DC bus} 
The voltages of converter AC bus $m$ and DC bus $n$ are coupled by the PWM’s amplitude modulation factor, which can be set according to the modulation mode. The voltage magnitude at converter AC
bus is upper bounded by \cite{Bahrami2017,Erickson2001}
\begin{equation}
\label{eq::voltcoup}
V_m \le \delta V_n,\;m\in\mathcal{N}_{ac}^{c},\;n\in\mathcal{N}_{dc}^{c},
\end{equation}
where $V_m$ denotes the voltage magnitude at converter AC bus $m\in\mathcal{N}_{ac}^{c}$, $V_n$ the voltage at converter DC bus $n\in\mathcal{N}_{dc}^{c}$. $\delta$ denotes the PWM’s amplitude modulation factor. 

\subsubsection{Active power coupling between the converter AC bus and DC bus} 
The coupling between the converter active power on AC side and DC side is described by 
\begin{equation}
\label{eq::actpcoup}
P_m+P_n+P^{loss}=0,\;m\in\mathcal{N}_{ac}^{c},\;n\in\mathcal{N}_{dc}^{c},
\end{equation}
where $P_m$ and $P_n$ denote the active power at converter AC bus $m\in\mathcal{N}_{ac}^{c}$ and DC bus $n\in\mathcal{N}_{dc}^{c}$, respectively. $P^{loss}$ denotes the converter power losses. 

The converter power losses can be approximated as a quadratic function of the phase current of the VSC valve as discussed in~\cite{Bahrami2017,Yang2018},
\begin{equation}
\label{eq::lossC}
P^{loss}=a_{1m} I_m^2+a_{2m} I_m +a_{3m},\;m\in\mathcal{N}_{ac}^{c}
\end{equation}
with 
\begin{equation}
\label{eq::lossI}
I_m =\sqrt{\frac{P_m^2+Q_m^2}{V_m^2}},\;m\in\mathcal{N}_{ac}^{c},
\end{equation}
where $I_m$ denotes the current magnitude at converter AC bus $m\in\mathcal{N}_{ac}^{c}$. $Q_m$ denotes the reactive power at converter AC bus $m\in\mathcal{N}_{ac}^{c}$. $a_{1m}$, $a_{2m}$, and $a_{3m}$ denote the conduction losses of the valves, the switching losses of valves and freewheeling diodes, and the no load losses of transformers and averaged auxiliary equipment losses, respectively. The values of the coefficients depend on the components and the power rating of the VSC station, and can be obtained using various approaches such as online identification or by aggregating the loss patterns of each component.

The converter losses contribute a relatively large percentage of the total system losses. The highly non-convex of converter losses adds a considerable computational burden for the AC/DC \acrshort{opf} problem in \acrshort{vsc-mtdc} meshed AC/DC grids.

\subsubsection{Apparent capacity limit of a converter} 
The power injected into a converter from AC side is calculated by
\begin{subequations}
\label{eq::Caceq}
\begin{align}
&P_m=V_m^2 g_m-V_mV_f \left(g_m\cos\theta_{mf}+b_m\sin\theta_{mf}\right)\\[0.12cm]
&Q_m=-V_m^2 b_m+V_mV_f \left(b_m\cos\theta_{mf}-g_m\sin\theta_{mf}\right)
\end{align}
\end{subequations}
for all $m\in\mathcal{N}_{ac}^{c}$ and $f\in\mathcal{N}_{ac}^{f}$. Here $g_m$ and $b_m$ denote the conductance and susceptance of phase reactor. 

The VSC station is principally constrained by the maximum current through VSC valves and the maximum DC voltage~\cite{2004Power}. The former determines the maximum VSC apparent power limit, and the latter defines the VSC reactive power output limit. In the present paper, VSC constraints at the converter AC bus are used because the VSC power exchange at the converter AC bus is set as control variables. The maximum current through the VSC valve has an upper limit as follows. 
\begin{equation}
\label{eq::Imax}
I_m\le \overline I_m,\;m\in\mathcal{N}_{ac}^{c},
\end{equation}
where $\overline I_m$ denotes the maximum current through the VSC valve at converter AC bus $m\in\mathcal{N}_{ac}^{c}$. 

Upon the substitution of \eqref{eq::lossI} into \eqref{eq::Imax}, the apparent capacity limit of the converter is expressed as
\begin{equation}
\label{eq::Smax}
P_m^2+Q_m^2\le \left(V_m \overline I_m\right)^2,\;m\in\mathcal{N}_{ac}^{c}.
\end{equation}

\subsubsection{Reactive power limit of the converter} 
The operation of the VSC station is constrained by the upper and lower limits of the converter reactive power output. In practical VSCs, the maximum reactive power that the converter absorbs is approximately proportional to the nominal value of its apparent power~\cite{Bahrami2017,Feng2014}. 
\begin{equation}
\label{eq::Qmax}
-\gamma S_m\le Q_m,\;m\in\mathcal{N}_{ac}^{c},
\end{equation}
where $\gamma$ denotes a positive constant and can be determined by the type of the converter, $S_m$ denotes the nominal value of the apparent power of the converter.

According to \cite{Bahrami2017,Yang2018}, since the susceptance $b_m$ of the phase reactor is normally much larger than its conductance $g_m$, the upper limit on the reactive power produced by the
converter is expressed as
\begin{equation}
\label{eq::Qmax3}
Q_m\le -b_m \overline V_m\left(\overline V_m- V_f\right),\; m\in\mathcal{N}_{ac}^{c},f\in\mathcal{N}_{ac}^{f}. 
\end{equation}
where $\overline V_m$ denotes the maximum voltage magnitude of converter AC bus $m\in\mathcal{N}_{ac}^{c}$, $V_f$ the voltage magnitude of AC filter bus $f\in\mathcal{N}_{ac}^{f}$.

\subsubsection{Nodal power balance of DC grid}
\begin{equation}
\label{eq::pbDC}
P_n=V_i\sum\limits_{\left({i,j}\right)\in {\mathcal{L}^{mtdc}}} {g_{ij}\left( {V_i - V_j} \right)},\; i \in {\mathcal{N}^{dc}},n\in\mathcal{N}_{dc}^{c},
\end{equation}
where $g_{ij}$ denotes the conductance of DC branch $(i,j)\in \mathcal{L}^{mtdc}$, $V_i$ denotes the nodal voltage of DC bus $i \in {\mathcal{N}^{dc}}$. 

\subsubsection{Branch flow limit of DC grid}
\begin{equation}
\label{eq::currDC}
-\overline P_{ij} \le {g_{ij}\left( {V_i- V_j} \right)V_i}  \le \overline P_{ij},\;\left( {i,j} \right) \in \mathcal{L}^{mtdc},
\end{equation}
where $\overline P_{ij}$ denotes the maximum capacity of DC branch $(i,j)\in \mathcal{L}^{mtdc}$.

\subsubsection{Nodal voltage limit of DC grid} 
The terminal voltage of one DC bus in DC grid fixed as a reference has to satisfy the limit
\begin{equation}
\label{eq::voltDC}
\underline V_i \le V_i \le \overline V_i,\;i \in {\mathcal{N}^{dc}}.
\end{equation}


\subsection{Centralized and Distributed Problem}
\subsubsection{Centralized Problem}
The full AC/DC \acrshort{opf} Problem of \acrshort{vsc-mtdc} meshed AC/DC grids is summarized as follows:
\begin{equation}
\label{eq::mainProb}
\begin{aligned}
&\text{minimize}\quad \eqref{eq::C}\\
\text{subject to}\;\;\;&
\eqref{eq::nodalAC}-\eqref{eq::Caceq},\eqref{eq::Smax}-\eqref{eq::voltDC}.
\end{aligned}
\end{equation}

\tb{Overall, this full AC/DC \acrshort{opf} model is a non-convex nonlinear optimization problem, which could have solutions based on local optimums. The interior point method is commonly used in solving large-scale nonlinear optimization problems. This paper focuses on the formulation of the distributed non-convex optimization rather than on the development of the advanced solution algorithms for the global optimum of the nonlinear optimization problem. For different (local optimal) solutions with the same objective function value, the total costs obtained based on local optimums could be the same; however, the generation dispatches and system losses could be different.}

For the proposed problem~\eqref{eq::mainProb}, it can only be directly solved in a centralized manner. However, the centralization will lead to the operation by a single central entity with complete knowledge and control of the entire meshed AC/DC grids. This centralization may create substantial regulatory and political issues because the local system operators have to give up their governance and control to the central entity. Thus, a distributed architecture is preferred. The distribution optimization repeatedly alternates solving a series small-scale subproblems and deploying a consensus step for coordination and hence, the iterative process is inevitable. However, this iterative process is worthwhile and essential as the distributed algorithm can not only preserve the information privacy and decision independence, but also comply with the philosophy of electricity market operations.


\subsubsection{Distributed Problem}
\begin{figure}[htbp!]
\centering
\includegraphics[width=\linewidth]{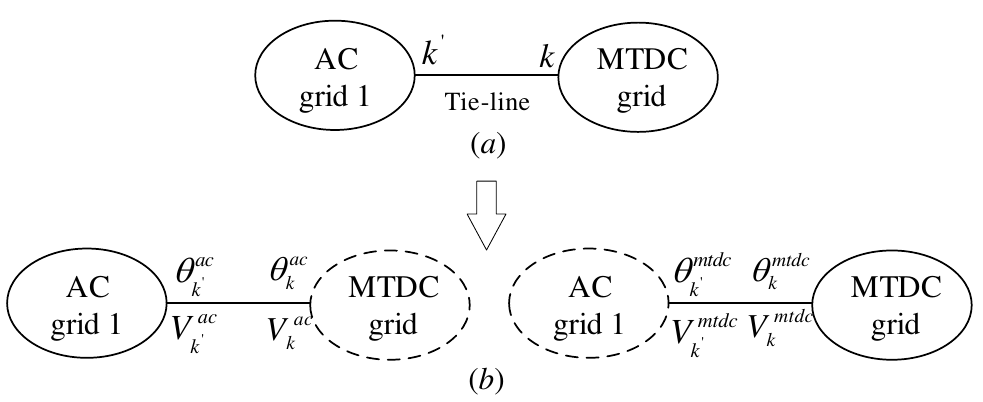}
\caption{Decomposition by sharing elements between neighboring regions}
\label{fig::Fg2}
\end{figure}

Similar to~\cite{muhlpfordt2021distributed}, we reformulate the full AC/DC \acrshort{opf} problem in affinely coupled separable form amenable to distributed optimization. Operating the \acrshort{vsc-mtdc} meshed AC/DC grids in a distributed way requires one to find the proper couplings between the regional AC grid and the MTDC grid. Take the coupling of AC grid 1 and MTDC grid as an example, shown in ~\autoref{fig::Fg2}(a), 
in which the system is divided into two parts that are connected through an AC tie-line between buses $k$ and $k'$. Buses $k$ and $k'$, and the linking AC tie-line are modeled together as a shared connection, labeled ($k$, $k'$), between these two parts. This shared connection is taken into account in the local problems of the AC grid and the MTDC grid as shown in~\autoref{fig::Fg2}(b). In the proposed distributed optimization algorithm, the nodal voltage magnitude and the phase angle on buses $k$ and $k'$ are coupled. Then, the consensus constraints are given by
\begin{subequations}
\label{eq::coup}
\begin{align}
V_k^{ac}&=V_k^{mtdc},\;\;V_{k'}^{ac}=V_{k'}^{mtdcc},\\
\theta_k^{ac}&=\theta_k^{mtdc},\;\;\theta_{k'}^{ac}=\theta_{k'}^{mtc},
\end{align}
\end{subequations}
which state that the complex voltage perceived by its connected zones should be identical.

\begin{remark}
The concept of sharing elements between neighboring regions allows one to compose the distributed \acrshort{opf} problem in a physically consistent manner: no additional modeling assumptions are introduced or required. If a solution to the distributed problem exists, then this will also be the solution to the respective centralized problem~\cite{6748974}\cite{muhlpfordt2021distributed}.
\end{remark}

For each AC grid $\ell\in\mathcal{R}$, the local objective is given by
\begin{equation}
\label{eq::Cj}
\begin{aligned}
f_\ell(x_\ell)=\sum_{i\in \mathcal{N}_\ell^{ac}}&\left\{c_{1i}\left( P_{i}^{G} \right)^{2}+c_{2i}P_{i}^{G}+c_{3i}\right\}\\
&\qquad\qquad\qquad +\eta\sum_{i\in \mathcal{N}_\ell^{ac}}\left(P_{i}^{G}-P_{i}^{D}\right),
\end{aligned}
\end{equation}
with stacked local variables $x_\ell$, while the objective of the \acrshort{mtdc} grid is set to be zero. This follows the fact that there is neither generator nor load in the \acrshort{mtdc} grid. Hence, generation cost $C_1= 0$ and losses term $C_2= 0$ for the \acrshort{mtdc} grid. 

Based on the description in this section, the full hybrid \acrshort{opf} problem \eqref{eq::mainProb} can be summarized into the standard affinely coupled distributed form
\begin{subequations}\label{eq::distOPT}
\begin{align}
\min_{x}\quad&f(x):=\sum_{\ell\in\mathcal{R}} f_\ell(x_\ell)\\\label{eq:affine1}
\textrm{s.t.}  \quad  &\sum_{\ell\in\mathcal{R}} A_\ell x_\ell =0\quad\;\;\mid\lambda\\\label{eq:affine2}
&h_\ell(x_\ell)\leq 0\qquad\quad\mid\kappa_\ell,\;\ell\in\mathcal{R}\\\label{eq:affine3}
&\underline x_\ell\leq x_\ell \leq \overline x_\ell\qquad\mid\gamma_\ell,\,\;\ell\in\mathcal{R}
\end{align}
\end{subequations}
Here, the affine constraints \eqref{eq:affine1} summarize the coupling~\eqref{eq::coup} between the AC grids and the MTDC grid. Constraints~\eqref{eq:affine2} and~\eqref{eq:affine3} collect all inequality constraints of each local model, where the box constraints~\eqref{eq:affine3} denote the bounds on the local generator active/reactive power and nodal voltage magnitudes.
In the following, we use notation
\[
\tilde{h}(x_\ell):=[h_\ell(x_\ell)^\top\;(x_\ell-\overline{x}_\ell)^\top\;(\underline{x}_\ell-x_\ell)^\top]^\top\leq 0,\;i\in\mathcal{R}
\] 
to stack~\eqref{eq:affine2} and~\eqref{eq:affine3}. Moreover, we use notation $x$ to concatenate $x_i$ for all $i\in\mathcal{N}$, and throughout the rest of this paper, we write down the Lagrangian multipliers right after the constraints such that $\lambda$, $\kappa_\ell$, $\gamma_\ell$ denote the dual variables (multipliers) of constraints~\eqref{eq:affine1},~\eqref{eq:affine2}, and~\eqref{eq:affine3}, respectively. 

\section{Distributed Optimization Algorithm}
\label{sec::singleArea}

This section presents \acrshort{admm} and  \acrshort{aladin} as algorithms to approach a distributed solution. Both algorithms share the same idea---update primal variables in an alternating fashion, whereas the major difference lays in the consensus step.

\subsection{ADMM}

The main idea of recalling \acrshort{admm} is to use it as a benchmark method reflecting the current state-of-the-art. In order to apply \acrshort{admm} for solving~\eqref{eq::mainProb}, we copy the variables $x_\ell$ and replace the consensus constraint~\eqref{eq:affine1} by  
\begin{equation}
x_\ell = z_\ell \quad \mid \xi_\ell,\;\;\ell\in\mathcal{R}\quad\text{and}\quad 0 = \sum_{i\in\mathcal{R}}A_\ell z_\ell
\end{equation}
with Lagrangian multipliers $\xi_\ell$, $\ell\in\mathcal{R}$.

Algorithm~\ref{alg::admm} outlines the main steps of~\acrshort{admm}. 
Step~\ref{alg::admm::s1a} solves $|\mathcal{R}|$ decoupled problems that are constructed according to the augmented Lagrangian. Step~\ref{alg::admm::s1b} updates the dual iterate $\xi_\ell$ based on the gradient ascend method~\cite{Boyd2011}. Notice that both Step~\ref{alg::admm::s1a} and~\ref{alg::admm::s1b} can be executed in parallel. Then, a practical terminal condition $\norm{\sum_{\ell\in\mathcal
R} A_\ell x_\ell^+}_\infty \leq \epsilon$
can be checked based on the local solution $x_\ell^+$.
Step~\ref{alg::admm::s2} deals a consensus QP whose solution can be worked out analytically, requiring one to collect $x_\ell^+$ and $\xi_\ell^+$. Once~\eqref{alg::admm::qp} is solved, the solution $z^+$ is broadcast to local agents, and the algorithm returns to Step~\ref{alg::admm::s1a}.
\begin{algorithm}[htbp!]
\caption{\acrshort{admm}\label{alg::admm}}
\textbf{Input:} $z$, $\xi$, $\rho>0$\\
\textbf{Repeat:}
\begin{enumerate}
\item update $x_\ell$ by solving all decoupled \acrshort{nlp} problems \label{alg::admm::s1a}
\begin{align}
	x_\ell^+=\underset{\tilde h_\ell(x_\ell)\leq0}{\mathrm{argmin}}\quad &f_\ell(x_\ell)+\xi_\ell^\top x_\ell+
	\frac{\rho}{2}\norm{x_\ell-z_\ell}^2_2
\end{align}
\item compute  \label{alg::admm::s1b}
$\xi_\ell^{+}=\xi_\ell+\rho (x_\ell-z_\ell)$, for all $\ell \in \mathcal{R}$.
\item update $z$ by solving the coupled averaging step~\label{alg::admm::s2}
\begin{subequations}\label{alg::admm::qp}
	hou\begin{align}
		z^+=\text{arg}\min_{z}\quad &\sum_{\ell\in\mathcal{R}} \left\{ \frac{\rho}{2} \norm{x_\ell^{+}-z_\ell}_2^2-\xi_\ell^{\top} z_\ell \right\}\quad\\
		\text{s.t.}\quad & \sum_{\ell\in\mathcal{R}} A_\ell z_\ell = 0
	\end{align}
\end{subequations}
\end{enumerate}
\end{algorithm}
\begin{remark}
By taking advantage of the topological structure, solving~\eqref{alg::admm::qp} in the consensus step~\ref{alg::admm::s2} of Algorithm~\ref{alg::admm} could only require neighbor-to-neighbor communications. This follows the fact that matrices $A_\ell$ denote the adjacency of the graph. More details refer to~\cite{6731604}.
\end{remark}

\subsection{ALADIN} 
\acrshort{aladin}, originally proposed in~\cite{Boris2016}, is developed for dealing with generic distributed nonlinear programming. Solving~\eqref{eq::distOPT} by using~\acrshort{aladin} is outlined in Algorithm~\ref{alg::aladin}. 
Similar to~\acrshort{admm}, Steps~\ref{alg::aladin::s1} and~\ref{alg::aladin::s2} of Algorithm~\ref{alg::aladin} are parallelizable. The local problems~\eqref{alg::aladin::nlp} are also formulated following the idea of the augmented Lagrangian. Here, one may adjust either the scaling matrices $\Sigma_\ell$ or penalty parameter $\rho$ during the iterations in order to improve performance. A practical strategy to update $\rho$ for distributed AC \acrshort{opf} can be found in~\cite{Engelmann2019}. Based on the local solutions $x_\ell$, the algorithm terminates if 
\begin{align*}
\norm{\sum_{\ell\in\mathcal{R}} A_\ell x_\ell }_\infty \leq \epsilon \;\text{and}\;
\max_\ell \norm{\Sigma_\ell(x_\ell-z_\ell)}_\infty \leq \epsilon. 
\end{align*}
holds. This condition implies that $x_\ell$ satisfies the first order optimality condition of~\eqref{eq::distOPT} 
\[
\left\|\nabla \left\{f_\ell(x_\ell)+\kappa^\top h_\ell(x_\ell) + \gamma_\ell^\top x_\ell\right\} + A_\ell^\top\lambda\right\| \leq \epsilon
\]
up to the user-specified numerical tolerance $\epsilon$. As discussed in~\cite{Boris2016}, the iterate $(x,\lambda,\kappa,\gamma)$  is a primal-dual KKT
point of Problem~\eqref{eq::distOPT}—up to the user specified numerical accuracy $\epsilon$. Step~\ref{alg::aladin::s2} evaluates the sensitivities at the local primal-dual solutions. \tb{Here, the Hessian approximation~\eqref{eq::Hessian} is required to be positive definite such that QP~\eqref{alg::aladin::qp} is convex and has unique solution. In practice, some numerical heuristics can be applied to make it be satisfied such as adding a small regularization. In our implementation, we adopt the heuristic implemented in open-source toolkit~\texttt{ACADO}~\cite{Houska2011}, which flips the negative eigenvalue of $H_\ell$.}

\begin{algorithm}[htbp!]
\caption{\acrshort{aladin}\label{alg::aladin}}
\textbf{Input}: $z$,\;$\lambda$,\;$\rho>0$,\;$\mu>0$ and scaling symmetric matrices $\Sigma_\ell\succ 0$\\
\textbf{Repeat:}
\begin{enumerate}
\item solve the following decoupled \acrshort{nlp}s for all \tb{$\ell\in\mathcal{R}$} \label{alg::aladin::s1}
\begin{subequations}
	\label{alg::aladin::nlp}
	\begin{align}
		\min_{x_\ell}\quad &f_\ell(x_\ell)+\lambda^\top A_\ell x_\ell+
		\frac{\rho}{2}\norm{x_\ell-z_\ell}^2_{\Sigma_\ell}\\
		\text{s.t.}\quad &h_\ell(x_\ell)\leq0 \qquad \mid\kappa_\ell,\\
		&\underline x_\ell\leq x_\ell \leq \overline x_\ell\quad\mid\gamma_\ell
	\end{align}
\end{subequations}
\item compute the Jacobian matrix $J_\ell$ based on the active set at the local solution $x_\ell$ by \label{alg::aladin::s2}
\begin{equation}
	[J_\ell]_i=\begin{cases}
		\partial\;[\tilde{h}_\ell(x_\ell)]_i
		&\text{if } [\tilde{h}_\ell(x_\ell)]_i=0\\[0.12cm]
		0&\text{otherwise}
	\end{cases}
	\label{eq::Jacobian}    
\end{equation}
with $[J_\ell]_i$ the i-th row of matrix $J_\ell$ and gradient $g_\ell=\nabla f_\ell(x_\ell)$, choose Hessian approximation
\begin{equation}\label{eq::Hessian}
	\tb{H_\ell\approx\nabla^2\left\{f_\ell(x_\ell)+\kappa_\ell^\top h_\ell(x_\ell)\right\}\succ 0. }
\end{equation}
\item  solve coupled \acrshort{qp} \label{alg::aladin::s3}
\end{enumerate}
\vspace{-0.2cm}
\begin{subequations}\label{alg::aladin::qp}
\begin{align}
	\min_{\Delta x,s}& \sum_{\ell\in\mathcal{R}}\left\{\frac{1}{2} \Delta x_\ell^\top\;H_\ell \;\Delta x_\ell + g_\ell^\top \Delta x_\ell\right\}
	+\lambda^\top\;s + \frac{\mu}{2} \norm{s}^2_2\\
	\textrm{s.t.}  & \quad  \sum_{\ell\in\mathcal{R}} A_\ell (x_\ell+ \Delta x) = b + s \quad \mid\lambda^{QP}\\\label{eq::Active}
	&  \quad J_\ell \Delta x_\ell = 0, \;\;\ell \in \mathcal{R}\;.
\end{align}
\end{subequations}
\begin{enumerate}
\item[4)] update primal and dual variable $z$, $\lambda$ by \label{alg::aladin::s4}
\begin{subequations}
	\begin{align}
		z^{+} &= x +\alpha_1(x-z) + \alpha_2 \Delta x, \\
		\lambda^{+} &= \lambda + \alpha_3 (\lambda^{QP}-\lambda), 
	\end{align}
\end{subequations}
where the line search scheme~\cite[Algorithm~3]{Boris2016} can be used to calculate the step sizes $\alpha_1$, $\alpha_2$ and $\alpha_3$.
\end{enumerate}
\end{algorithm} 
Notice that the size of Jacobian matrices $J_\ell$ is changed during the iterations as the active set~\footnote{The active set of~\eqref{alg::aladin::nlp} is defined by $\{i\in\mathbb{N}_+\,\mid\,[\tilde h_\ell(x_\ell)]_i=0,\forall i\}$.}might be changed. In practice, if the quasi-Newton Hessian approximation such as Broyden–Fletcher–Goldfarb–Shanno (\textsc{bfgs}) is used to compute $H_\ell$, the communication cost can be significantly reduced~\cite{Engelmann2019}. The main difference between~\acrshort{admm} and~\acrshort{aladin} is the consensus step. While both approaches require one to solve a \acrshort{qp} problem, different from~\eqref{alg::admm::qp}, Problem~\eqref{alg::aladin::qp} is equivalent to an inexact Newton step for solving~\eqref{eq::distOPT} with only considering the active constraints at the local solution $x_\ell$. This is crucial for the convergence improvement of~\acrshort{aladin}. Notice that introducing the slack variables $s$ in~\eqref{alg::aladin::qp} guarantees that Problem~\eqref{alg::aladin::qp} is always feasible no matter if the original problem~\eqref{eq::distOPT} is feasible or not. \tb{In practice, solving the equality constrained QP~\eqref{alg::aladin::qp} only needs a basic linear algebraic routine such as~\texttt{Lapack} as it is equivalent to solve the resulting KKT system based linear equations. Its analytical solution can be found in Appendix. One can see that if the full step is applied at Step~4, only the dual update~\eqref{eq::dualSol} requires communication while it is not necessary to communicate all sensitivities. Thus, if there exists a central entity or any local agent could perform as a central coordinator, Step 3 can be efficiently implemented. Moreover, if only neighbor-to-neighbor communication is allowed,~\cite{Engelmann2020} has proposed bi-level distributed variants of \acrshort{aladin}, in which three methods, including Schur complement method, decentralized ADMM method and decentralized conjugate gradient method, was adopted to deal with~\eqref{alg::aladin::qp} in a decentralized manner.}

The overall frameworks of Algorithms~\ref{alg::admm} and~\ref{alg::aladin} are similar. However, in general, applying~\acrshort{admm} to nonconvex optimization~\eqref{eq::distOPT} does not have any theoretical guarantees. In contrast to this,~\acrshort{aladin} has a local convergence guarantee while the globalization presented in~\cite[Algorithm~3]{Boris2016} can guarantee the global convergence with doing line search in Step~\ref{alg::aladin::s4} of Algorithm~\ref{alg::aladin}. In this paper, we focus on the local convergence properties of Algorithm~\ref{alg::aladin} as in practice, according to the physical model of the power grids, a good initial guess of the primal-dual iterates can be obtained. In order to introduce the local convergence results of Algorithm~\ref{alg::aladin}, we need the following definition.
\begin{definition}
~A KKT point for generic constrained optimization is called regular~\cite{nocedal2006numerical} if the linear independence constraint qualification (LICQ), strict complementarity
conditions (SCC), as well as the second order sufficient condition (SOSC) are satisfied. 
\end{definition}
Now, we summarize the local convergence property of Algorithm~\ref{alg::aladin} as follows:
\begin{theorem}
\label{thm::convergence}
Let the KKT point $(z^*,\lambda^*,\kappa^*,\gamma^*)$ of Problem~\eqref{eq::distOPT} be regular such that following the SOSC it is a local minimizer. And let the penalty parameter $\rho<\infty$ in~Algorithm~\ref{alg::aladin} be sufficiently large with $\nabla^2\left\{f_\ell(x_\ell)+\kappa_\ell^\top h_\ell(x_\ell)\right\}+\rho \Sigma_\ell\succ 0$. Moreover, let matrices
\begin{equation}\label{eq::HessianCond}
H_\ell = \nabla^2\left\{f_\ell(x_\ell)+\kappa_\ell^\top h_\ell(x_\ell)\right\} + \mathcal{O}(\|x_\ell-z_\ell\|) 
\end{equation}
holds for all $\ell\in\mathcal{R}$, the iterates $x_\ell$ of Algorithm~\ref{alg::aladin} converge locally with quadratic rate if full step size is appled in Step 4), i.e., $\alpha_1=\alpha_2=\alpha_3=1$.
\end{theorem}
Here, the local convergence means that the initial guess of primal-dual iterates are located in a small neighborhood at the local minimizer $(x^*=z^*,\lambda^*)$. The proof of Theorem~\ref{thm::convergence} can be established in two steps. First, according to the assumptions of regularity and $\rho$, \tb{the local minimizer of subproblems~\eqref{alg::aladin::nlp}, $x_\ell$ is parametric with $(z,\lambda)$ and the solution maps are Lipschitz continuous, i.e.,
there exists constants $\chi_1,\chi_2>0$ such that 
\begin{equation}\label{eq::Lipschitz}
\left\|x-z^*\right\|\leq  \chi_1\left\|z-z^*\right\| + \chi_2\left\|\lambda-\lambda^*\right\|.
\end{equation}
This result was formally stated in~\cite[Lemma~3]{Boris2016} and the proof can be established by applying the implicit function theorem~\cite[Appendix~2, Page 630]{nocedal2006numerical}, which refer to~\cite[Lemma~1]{Houska2021} for more details.
The second step follows the fact that the active sets are not changed locally based on the assumptions. Then, the standard analysis of Newton's method~\cite[Chapter 3.3]{nocedal2006numerical} gives 
\[
\begin{aligned}
\left\|\begin{bmatrix}
z^+-z^*\\ \lambda^+-\lambda^*
\end{bmatrix}\right\|\leq& \left\|H-\nabla^2f(x)+\kappa^\top h_\ell\right\|\cdot\mathcal{O}\left(\left\|x-z^*\right\|\right)\\ &+\mathcal{O}\left(\left\|x-z^*\right\|^2\right)	
\end{aligned}
\]
Based on condition~\eqref{eq::HessianCond}, we have that there exists a constant $\sigma>0$ such that the local quadratic contraction
\begin{equation}\label{eq::contraction}
\left\|z^+-z^*\right\|\leq  \sigma\left\|x-z^*\right\|^2 \;\;\text{and}\;\;
\left\|\lambda^+-\lambda^*\right\|\leq  \sigma\left\|x-z^*\right\|^2
\end{equation}
can be established. Combing~\eqref{eq::Lipschitz} and~\eqref{eq::contraction} yields the result. A more detailed proof of Theorem~\ref{thm::convergence} can be found in~\cite[Section 7]{Boris2016}. }
\begin{remark}
If we use BFGS Hessian approximation, a local superlinear convergence rate is achieved instead~\cite{jiang2021block}. The analysis is close to Theorem~\ref{thm::convergence} with slight difference in the second step of the proof.  
\end{remark}
\begin{remark}
\tb{In practice, only a numerical solution of decoupled problem~\eqref{alg::aladin::nlp} is achievable while convergence result in Theorem~\ref{thm::convergence} can be still established by assuming that the local solutions of~\eqref{alg::aladin::nlp} is bounded by $\left\|\bar{x} - x\right\|\leq \zeta \left\|x-z\right\|$ with $\bar x$ the approximate solution given by the numerical solver.   The proof refers to~\cite[Theorem~2]{Engelmann2019}}
\end{remark}

\section{Numerical Results}
\label{sec::result}

\begin{figure*}[htbp!]
\centering
\includegraphics[width=0.90\textwidth]{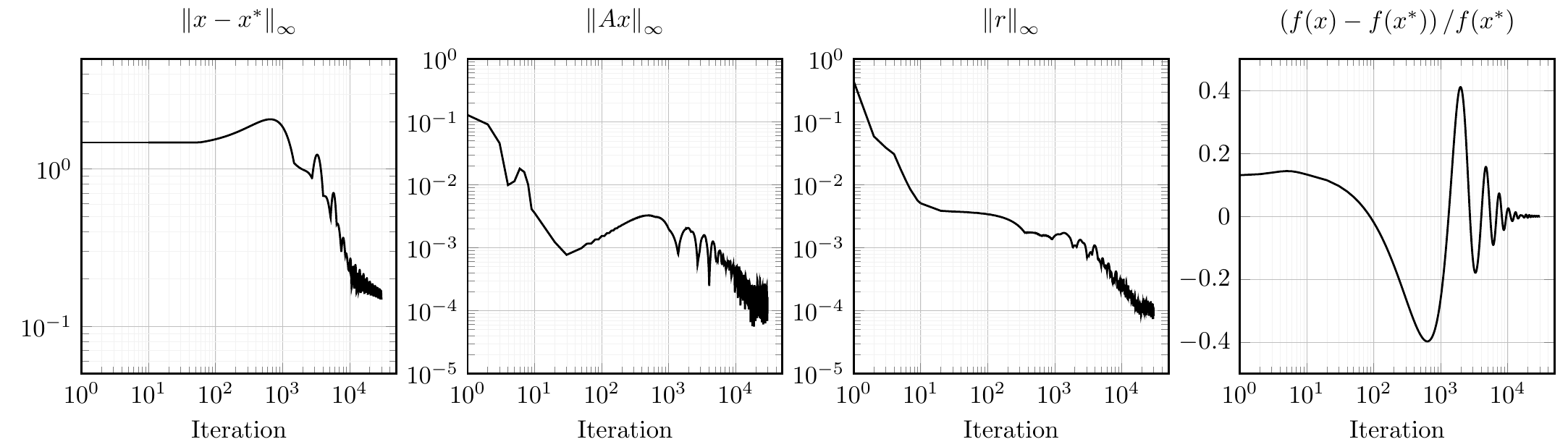}
\caption{Convergence behavior of \acrshort{admm} for Case 1}\label{fig::convergence::admm}
\end{figure*}

\begin{figure*}[htbp!]
\centering
\includegraphics[width=0.90\textwidth]{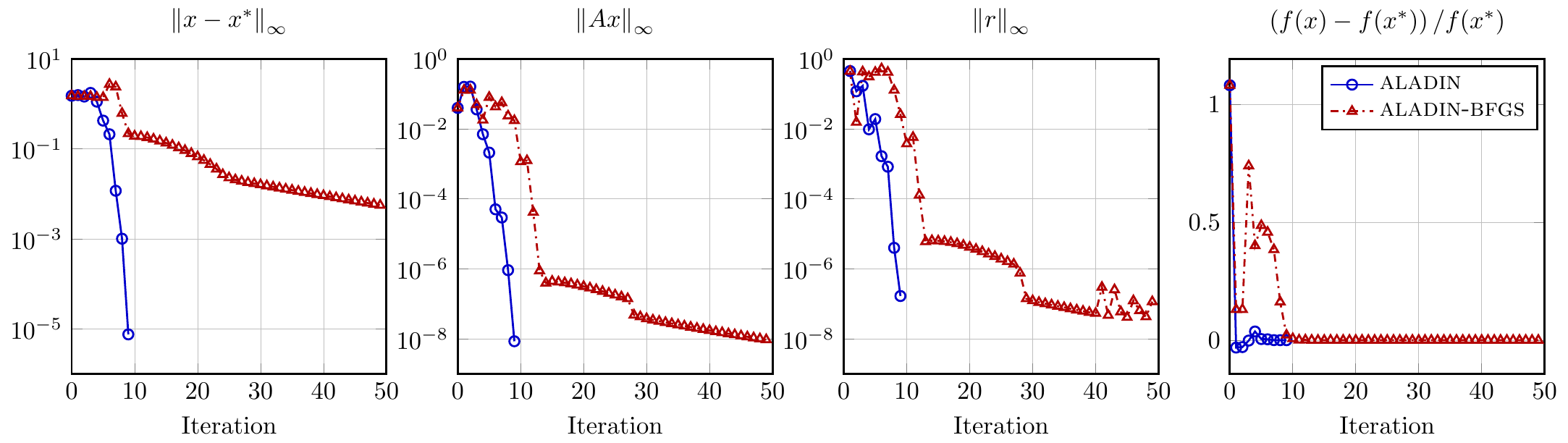}
\caption{Convergence behavior of \acrshort{aladin} for Case 1}\label{fig::convergence::aladin::case1}
\end{figure*}

\begin{figure*}[htp]
\centering
\includegraphics[width=0.90\textwidth]{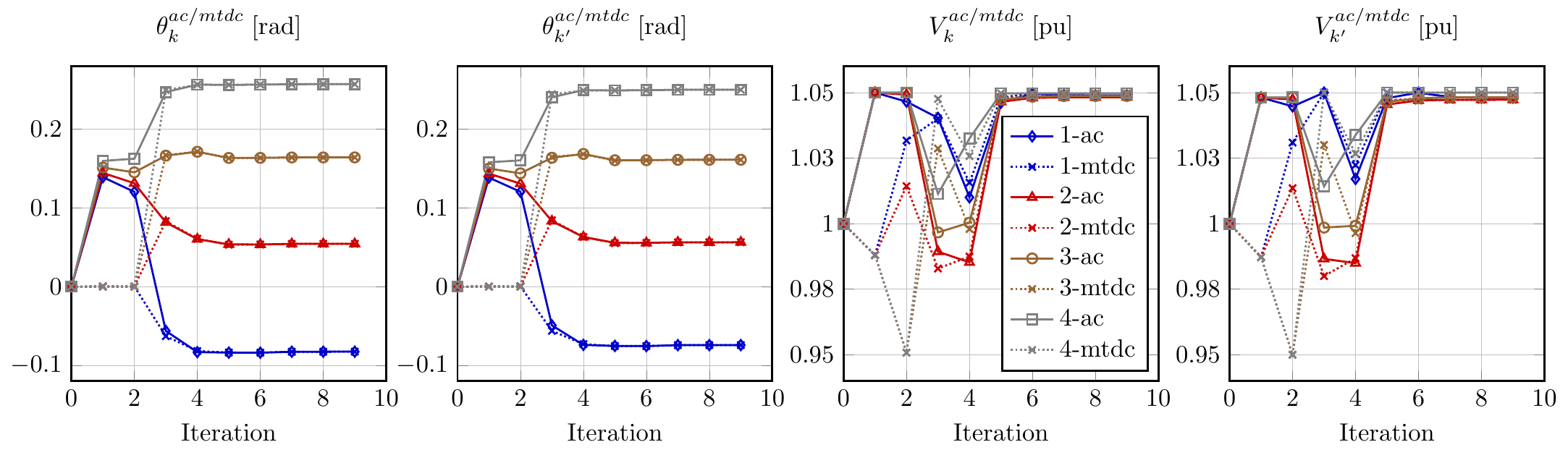}
\caption{Coupling variables convergence towards optimal value of \acrshort{aladin} with exact Hessian for Case 1}\label{fig::coupling::aladin}
\end{figure*}

\begin{figure*}[htp]
\centering
\includegraphics[width=0.90\textwidth]{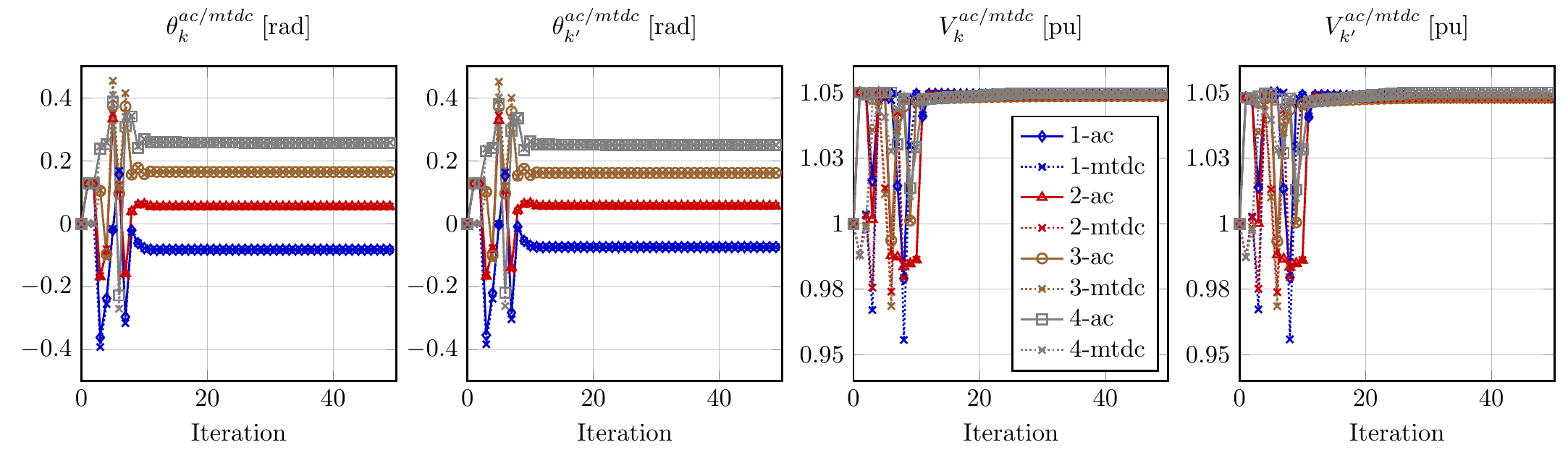}
\caption{Coupling variables convergence towards optimal value of \acrshort{aladin} with inexact Hessian for Case 1}\label{fig::coupling::bfgs}
\end{figure*}

\begin{figure*}[htbp!]
\centering
\includegraphics[width=0.90\textwidth]{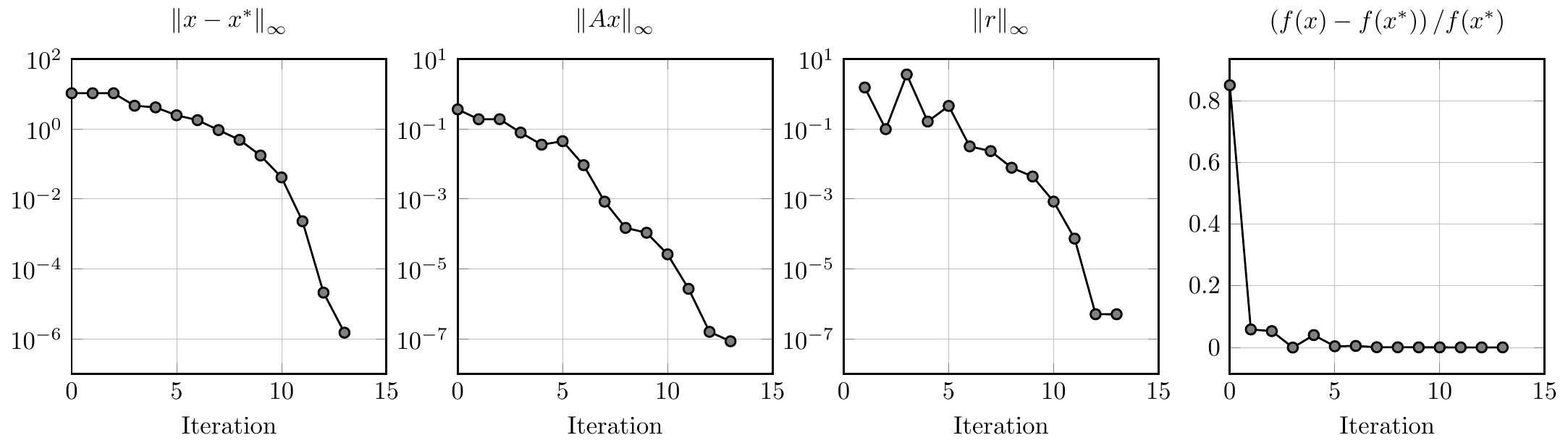}
\caption{Convergence behavior of \acrshort{aladin} using exact Hessians for Case 2}\label{fig::convergence::case2} 
\end{figure*}

This section illustrates the performance of the proposed distributed algorithms with comparison to \acrshort{admm}. The model is merged by a four-terminal \acrshort{vsc}-\acrshort{hvdc} network with four AC grids from \matpower~\cite{zimmerman2010matpower}, which are modified using open-source \acrshort{rapidpf} toolbox\footnote{
The code is available on \url{https://github.com/KIT-IAI/rapidPF}}, proposed in~\cite{muhlpfordt2021distributed}.

In Case 1 a \acrshort{vsc}-\acrshort{hvdc} network is connected with four IEEE 9-bus AC grids as shown in~\autoref{Fg::vsc-mtdc} while a \acrshort{vsc}-\acrshort{hvdc} network is connected with four IEEE 118-bus AC grids in Case 2. Bus 2 and Bus 8 are the connecting bus for 9-bus AC grid and 118-bus AC grid respectively. For both Cases, base power and base voltage is set to 100MVA and 345kV; voltage magnitude of both AC and DC grid are limited to [0.95, 1.05] \acrshort{pu}; the master VSC station connecting to AC grid 3 keeps 1 p.u. to provide the constant DC voltage reference; the parameters for \acrshort{vsc} stations and DC branches in the grid’s \acrfull{pu} are given in ~\autoref{TB::vsc-station} and \ref{TB::DC-grid} to differentiate the four AC grids in both Cases, the generator cost coefficients and load of AC regions 1 and 2 are set to be smaller than regions 3 and 4, shown in \autoref{TB::AC-setting}. In this way, we force the power exports from AC regions 1 and 2 to regions 3 and 4. The scaling coefficient for system loss term is set to $\eta=10$ to jointly minimize the total generation cost and system losses. \tbb{For reasons of a fair comparison, our implementation initializes the primal variables $x$ for each AC block with 1 p.u for voltage magnitudes while 0 for all other values. This flat starting strategy is standard, which was used in literature, for example,~\cite{Engelmann2019,muhlpfordt2021distributed}, and also chosen as the default strategy in~\texttt{Matpower}~\cite{zimmerman2010matpower}. For the \acrshort{mtdc} block, we applied a similar strategy that intuitively initialized voltage magnitudes by $1$ p.u, and the injected power by $0$.} The initial dual variable $\lambda$ are set to zero (flat start).

Applying \acrshort{admm} and \acrshort{aladin} to solve problems requires selecting tuning parameters $\rho$ and $\mu$. To enable a fair comparison, these parameters are determined by parameter sweeps aiming for fast convergence, shown in Table \ref{TB::tuning-parameter}. To obtain a similar scaling, the weighting matrices are chosen such that each diagonal entry is inversely proportional to its corresponding decision variable range, as suggested in~\cite{Engelmann2019}. 
\begin{table}[htbp!]
\caption{VSC station parameters} 
\label{TB::vsc-station}
\centering
\scriptsize
\renewcommand{\arraystretch}{1.4}
\begin{tabular}{cccccc}                           \hline
Term   & Value & Term   & Value & Term     & Value  \\\hline
$a_{1m}$ & 0.011 p.u.  &$R_m$  &  0.00025 p.u.  & $R_f$  &  0.0005 p.u.       \\
$a_{2m}$ & 0.003 p.u.  &$X_m$  &  0.04 p.u.     & $X_f$ &   0.0125 p.u.    \\
$a_{3m}$& 0.0043 p.u.  &$B_f$  & 0.2 p.u.       &$S_m$  &   11 p.u.       \\
$\delta$&  1.05 &$\gamma$ &  0.5 & $\overline V_m$   & 1.05 p.u.\\
$\overline I_m$   &  11 p.u.\\      \hline
\end{tabular}
\end{table}

\begin{table}[htbp!]
\caption{DC branch parameters} 
\label{TB::DC-grid}
\centering
\scriptsize
\renewcommand{\arraystretch}{1.4}
\begin{tabular}{cccccc}
\hline
\multirow{2}{*}{Line No.} &\multirow{2}{*}{From bus}  &\multirow{2}{*}{To bus}&  \multirow{2}{*}{Resistance [pu]} & \multicolumn{2}{c}{Flow limit [pu]}\\
&          &        &              & Case 1  & Case 2        \\\hline		
1        & 1        & 2      &  0.00042     & 1.5     & 8            \\
2        & 1        & 3      &  0.00174     & 1.5     & 8            \\
3        & 2        & 4      &  0.00175     & 1.5     & 8          \\
4        & 3        & 4      &  0.00159     & 1.5     & 8           \\\hline          
\end{tabular}
\end{table}


\begin{table}[htbp!]
\caption{Load and generator cost coefficients ratio setting} \label{TB::AC-setting}
\centering
\scriptsize
\renewcommand{\arraystretch}{1.4}
\begin{tabular}{ccccc}                           \hline
&AC grid No. & $P^D$ [MW] & $Q^D$ [MVar] & Cost Coefficients ratio\\ \hline
\multirow{4}{*}{Case 1}& 1 & 315  & 115   & 1.0\\
& 2 & 318  & 116   & 1.3\\
& 3 & 321  & 117   & 1.7\\
& 4 & 324  & 118   & 2.2\\\hline
\multirow{4}{*}{Case 2}& 1 & 4242 & 1438  & 1.0\\
& 2 & 4666 & 1582  & 1.3\\
& 3 & 5090 & 1726  & 1.7\\
& 4 & 5515 & 1869  & 2.2\\\hline
\end{tabular}
\end{table}

\begin{table}[htbp!]
\caption{Algorithm parameters of \acrshort{admm} and \acrshort{aladin}} \label{TB::tuning-parameter}
\centering
\scriptsize
\renewcommand{\arraystretch}{1.4}
\begin{tabular}{ccccccc}
\hline
\multirow{2}{*}{Parameters} & \multicolumn{2}{c}{ADMM} & \multicolumn{2}{c}{ALADIN-BFGS} & \multicolumn{2}{c}{ALADIN}\\
& Case 1  & Case 2 & Case 1 & Case 2 & Case 1 & Case 2 \\\hline
$\rho$ & $10^4$  & - & $10^4$ & - & $10^2$ & $10^2$ \\
$\mu$  & - & - & $10^3$ & - & $10^3$ & $10^3$       \\\hline
\end{tabular}
\end{table}

The framework is built on \matlab-R2021a and the AC grid modeling follows AC-\acrshort{opf} model of \matpower. The case studies are carried out on a standard desktop computer with \texttt{Intel\textsuperscript{\textregistered} 
i5-6600K CPU @ 3.50GHz} and 16.0 \textsc{GB} installed \textsc{ram}. \casadi toolbox~\cite{andersson2019casadi} is used in \matlab and \ipopt~\cite{wachter2006implementation} is used as the solver for both the decoupled \acrshort{nlp}s and the coupled \acrshort{qp}. The computation time is posted in \autoref{TB::numerical-result}.

\begin{table*}[ht]
\caption{Comparisons of different algorithms for Cases 1 and 2} \label{TB::numerical-result}
\centering
\scriptsize
\renewcommand{\arraystretch}{1.5}
\begin{tabular}{ccccccccc}\hline
\multirow{2}{*}{Case}   & \multirow{2}{*}{Algorithm} & \multirow{2}{*}{Iterations} & \multirow{2}{*}{Time [s]}  & \multirow{2}{*}{$\norm{x-x^*}_\infty$}& \multicolumn{2}{c}{Generation Cost $C_1$}&  \multicolumn{2}{c}{System Losses $C_2$}\\
&                           &       &        &     & Cost  [$\times10^2 \$$] & Solution gap & Losses [MW] & Solution gap\\\hline
\multirow{4}{*}{1} & Centralized               & -  & $0.181$& -   & $249.354$ & - & $23.484$ & -\\
& \acrshort{admm}           &$\geq3\times 10^4$&$\geq3000$& $0.151$  & $249.141$& $8.56\times10^{-4}$ &$23.505$ & $8.95\times10^{-4}$\\
& \acrshort{aladin}-\acrshort{bfgs} & $49$  & $ 0.700$&$5.62\times 10^{-3}$&$249.354$&$2.03\times10^{-8}$&$23.484$& $1.77\times 10^{-6}$\\
& \acrshort{aladin}  & $ 9$  & $0.142$&$7.52\times10^{-6}$&$249.354$&$7.94\times10^{-7}$&$23.484$& $1.978\times 10^{-7}$\\\hline
\multirow{4}{*}{2} & Centralized               & -  &$1.465$ & -   & $9547.316$ &   -       & $515.985$ &  -\\
& \acrshort{admm}         & \multicolumn{7}{c}{Did not converge}\\
& \acrshort{aladin}-\acrshort{bfgs}         &  \multicolumn{7}{c}{Did not converge}\\
& \acrshort{aladin}  & $13$  &$0.878$ &$1.50\times 10^{-6}$&$9547.316$&$6.63\times10^{-9}$&$515.985$&$7.59\times10^{-9}$ \\\hline
\end{tabular}
\end{table*}
Similar to \cite{Engelmann2019}, we use the following quantities in order to depict the algorithm convergence behavior. The reference solution $x^*$ is obtained by solving problem \eqref{eq::mainProb} with \ipopt centrally.
\begin{enumerate}
\item The deviation of optimization variables from the optimal value $\norm{x-x^*}_\infty$.
\item The consensus violation $\norm{A x}_\infty = \norm{\sum_{\ell\in\mathcal{R}}A_\ell x_\ell}_\infty$.
\item The dual residual $\norm{r}_\infty$ with $r = \sum_{\ell \in\mathcal{R}}\left\{A_\ell(x_\ell-z_\ell)\right\}$ for \acrshort{admm}, $r = \sum_{\ell \in\mathcal{R}}\left\{\Sigma_{\ell}(x_\ell-z_\ell)\right\}$ for \acrshort{aladin}. 
\item The solution gap between a specific distributed algorithm and the centralized algorithm is calculated as $\frac{f(x)-f(x^*)}{f(x^*)}$.
\end{enumerate}
\begin{remark}
When the problem is feasible, it means that the \acrshort{vsc-mtdc} Meshed AC/DC Grid has at least one safe operating point, in which neither the transmission lines nor the converters are overloaded, etc. Meanwhile, the correct solution represents the operating point with the best economic efficiency for the grid.
\end{remark}

\subsection{Four-Terminal VSC-HVDC with Four IEEE 9-bus AC Grids}

\subsubsection{ADMM vs. ALADIN}

In general, \acrshort{admm} has no convergence guarantee for and AC-\acrshort{opf} problem, especially when the AC branch flow limits are implemented as the nonlinear inequality constraint\cite{Engelmann2019}. The mathematical model of \acrshort{vsc-mtdc} meshed AC/DC grids is much more complicated and the convergence of \acrshort{admm} is highly sensitive to the tuning parameters. \autoref{fig::convergence::admm} shows the convergence behaviors of \acrshort{admm}, with tuning parameter $\rho=10^{4}$. It is observed that \acrshort{admm} requires thousands of iterations to converge to modest accuracy, i.e., $\epsilon=10^{-4}$. To approach higher accuracy, primal variables $y$ start damping near the optimizer. Nevertheless, the solution gap by applying \acrshort{admm} is acceptable for practical use. 

In contrast to \acrshort{admm}, \acrshort{aladin} obtains curvature information of all subproblems in the central coordinator and has the ability to deal with the non-convex constraints. \autoref{fig::convergence::aladin::case1} presents convergence behavior of \acrshort{aladin} with two different methods for computing the Hessian matrix. It is obvious that \acrshort{aladin} converges significantly faster than \acrshort{admm}. Unlike \acrshort{admm}, the solution gap of \acrshort{aladin} is fairly small compared with the centralized method, as shown in Table \ref{TB::numerical-result}. The number of iterations and the computation time in \acrshort{aladin} is much less than \acrshort{admm}. The computational costs per iteration are only slightly higher for \acrshort{aladin}, which leads to a strong computation time decrease while much higher accuracy in terms of the optimal gap and dual residual are obtained. In summary, despite higher communication effort, \acrshort{aladin} surpasses \acrshort{admm} in all perspective for solving the AC/DC \acrshort{opf} problem of \acrshort{vsc-mtdc} meshed AC/DC grids. 

\subsubsection{Exact Hessian vs. BFGS}
To reduce the per-step communication effort between decoupled \acrshort{nlp}s and the central coordinator,  the \acrshort{bfgs} method is implemented to avoid communicating the full Hessian matrix.
\tb{Considering the theoretically worst case, i.e., all matrices and vectors required to communicate are dense without zero elements, Algorithm~\ref{alg::aladin} with exact Hessian requires to communicate 
\[
\sum_{\ell\in\mathcal{R}}\underbrace{n_\ell}_{g_\ell} +\underbrace{\frac{n_\ell(1+n_\ell)}{2}}_{\text{symmetric}\;H_i}+\underbrace{n_\ell \sum_{j\in\mathcal{A}_\ell}m_i}_{J_\ell}
\] 
float numbers in forward communication while using BFGS update can reduce this number as $\sum_{\ell\in\mathcal{R}}\left(3n_\ell + n_\ell \sum_{i\in\mathcal{A}_\ell}m_i \right)$. Here, $n_\ell$ defines the dimension of $x_\ell$ and set $\mathcal{A}_\ell$ collects the index of the active inequality constraints $[\tilde h_\ell(x_\ell)]_i=0$ at local solution $x_\ell$. } \autoref{fig::convergence::aladin::case1} compares the convergence behaviors of using exact Hessian and using \acrshort{bfgs} to approximate Hessian for solving the AC/DC \acrshort{opf} problem. For both options, \acrshort{aladin} can converge to a solution with modest accuracy, i.e., $\epsilon=10^{-4}$ for all termination criteria, within a dozen iterations. Furthermore, \acrshort{aladin} using inexact Hessians needs just slightly more iterations compared with \acrshort{aladin} using exact Hessians. One can observe that the convergence rate for \acrshort{aladin} using inexact Hessians is faster than linear. Nonetheless, it needs several dozens more iterations to converge when a high accurate solution is required, while \acrshort{aladin} with exact Hessian needs only several iterations. 

In the perspective of coupling variables, as illustrated in \autoref{fig::coupling::aladin} and \ref{fig::coupling::bfgs}, the nodal voltage magnitude and phase angle of boundary bus and fictitious bus among the local AC grids and MTDC grid converge to the same operating point after $9$ iterations (exact Hessians) and $49$ iterations (inexact Hessians), respectively.  

\subsection{Four-Terminal VSC-HVDC with Four IEEE 118-bus AC Grids}
The performance of the proposed algorithm is further analyzed on a larger \acrshort{vsc-mtdc} meshed AC/DC grid in Case 2. As summarized in \autoref{TB::numerical-result}, the \acrshort{admm} and \acrshort{aladin}-\acrshort{bfgs} fail to converge for solving this larger system no matter how the step size $\rho$ is adapted. 
The algorithm convergence behavior for Case 2 using \acrshort{aladin} with exact Hessians is depicted in \autoref{fig::convergence::case2}. 
Similar to Case 1, the solution gap of \acrshort{aladin} is fairly small compared with the centralized method. 

The proposed \acrshort{aladin} algorithm with exact Hessian takes 0.142 and 0.878 seconds for the two cases, respectively, and the computation time of \acrshort{aladin} is slightly faster than the centralized approach that takes 0.181 and 1.465 seconds, respectively. In conclusion, \acrshort{aladin} with exact Hessian outperform both centralized approach and \acrshort{admm} method for the two cases. This indicates that \acrshort{aladin} has more potential in practical distributed implementation.

\section{Conclusions}
\label{sec::conc}
To fully coordinate the AC grids and the MTDC grid to mutually benefit multiple regions, the present paper tailors the recently proposed novel distributed \acrshort{aladin} algorithm to solve the non-convex AC/DC OPF problem for the VSC-MTDC meshed AC/DC grids. By operating the AC grids and \acrshort{mtdc} grid in a fully distributed way, the information privacy and decision independency of local systems can be achieved, which is under the operating philosophy of the electricity market and the hierarchical and partitioned operating mode in China. In contrast to the commonly used \acrshort{admm}, \acrshort{aladin} has theoretical convergence guarantee and is able to outperform \acrshort{admm} in all perspective for the non-convex AC/DC \acrshort{opf} problem. 

\appendix
In order to work out the analytical solution of~\eqref{alg::aladin::qp}, we write down the KKT system of~\eqref{alg::aladin::qp} as follows
\[
\begin{bmatrix}
H&&A^\top & J^\top\\
&\mu I &-I & \\
A& -I & & \\
J & &&
\end{bmatrix}
\begin{bmatrix}
\Delta x \\ s \\ \lambda^{QP} \\ \kappa^{QP}
\end{bmatrix}
=\begin{bmatrix}
-g\\-\lambda\\b-Ax\\0
\end{bmatrix}
\]
with $H=\text{diag}\{H_\ell\}_{\ell\in\mathcal{R}}$, $J=\text{diag}\{J_\ell\}_{\ell\in\mathcal{R}}$, and $\Delta x$, $g$ and $\kappa^{QP} $ stacks $\Delta x_\ell$, $g_\ell$ and $\kappa^{QP}_\ell$ for all $\ell\in\mathcal{R}$ as column vectors. Here $\kappa^{QP}_\ell$ defines the Lagrangian multipliers of constraint~\eqref{eq::Active}.  By using the Schur complement, we can have the dual solution $\lambda^{QP}$ in a form
\begin{align}\label{eq::dualSol}
\lambda^{QP} &=  \left(
\sum_{\ell\in\mathcal{R}}
\begin{bmatrix}
A_\ell^\top\\0
\end{bmatrix}^\top
\begin{bmatrix}
H_\ell & J_\ell^\top\\J_\ell
\end{bmatrix}^{-1}
\begin{bmatrix}
A_\ell^\top\\0
\end{bmatrix}+\mu^{-1} I 
\right)^{-1}\\\notag
&\cdot\left(\mu^{-1}\lambda + 
\sum_{\ell\in\mathcal{R}}
\begin{bmatrix}
A_\ell^\top\\0
\end{bmatrix}^\top
\left(\begin{bmatrix}
x_i \\ 0
\end{bmatrix}
- 
\begin{bmatrix}
H_\ell & J_\ell^\top\\J_\ell
\end{bmatrix}^{-1}
\begin{bmatrix}
g_\ell\\0
\end{bmatrix}\right)
\right)
\end{align}
and the decoupled solutions for all $\ell\in\mathcal{R}$, 
\begin{equation}
\label{eq::deSol}
\begin{bmatrix}
\Delta x_\ell \\ \kappa_\ell^{QP}
\end{bmatrix}
=\begin{bmatrix}
H_\ell & J_\ell^\top\\J_\ell&
\end{bmatrix}^{-1}\left(
\begin{bmatrix}
-g_\ell\\0
\end{bmatrix}-
\begin{bmatrix}
A_\ell^\top\\0
\end{bmatrix}
\lambda^{QP}
\right).
\end{equation}

\bibliographystyle{ieeetr}
\bibliography{main}

\end{document}